\documentclass[final,1p,times]{elsarticle}
\usepackage{epsfig,amssymb}
\usepackage{amsfonts,amsmath,amsthm}

\def\eps{\epsilon}

\def\E{{\bf E}}
\def\H{{\bf H}}
\newcommand{\be}{\begin{equation}}
\newcommand{\ee}{\end{equation}}
\newcommand{\bea}{\begin{eqnarray}}
\newcommand{\eea}{\end{eqnarray}}

\newcommand{\x}{{\mathbf x}}
\newcommand{\y}{{\mathbf y}}

\newcommand{\xy}{\mathbf x}

\newcommand{\Sc}{{\cal S}}
\newcommand{\Dc}{{\cal D}}
\newcommand{\Nc}{{\cal N}}
\newcommand{\Tc}{{\cal T}}
\newcommand{\SQP}{{\cal S}_\tbox{QP}}
\newcommand{\DQP}{{\cal D}_\tbox{QP}}
\newcommand{\NQP}{{\cal N}_\tbox{QP}}
\newcommand{\TQP}{{\cal T}_\tbox{QP}}
\newcommand{\gqp}{{G}_\tbox{QP}}
\newcommand{\gqpnear}{{G}^\tbox{near}_\tbox{QP}}
\newcommand{\gqpfar}{{G}^\tbox{far}_\tbox{QP}}

\newcommand{\tbox}[1]{{\mbox{\tiny #1}}}

\newcommand{\Amat}{\mathsf{A}}
\newcommand{\Dmat}{\mathsf{D}}
\newcommand{\Lmat}{\mathsf{L}}
\newcommand{\Smat}{\mathsf{S}}
\newcommand{\Rmat}{\mathsf{R}}
\newcommand{\Imat}{\mathsf{I}}
\newcommand{\xmat}{\mathsf{x}}
\newcommand{\ymat}{\mathsf{y}}
\newcommand{\zmat}{\mathsf{z}}
\newcommand{\bmat}{\mathsf{b}}
\newcommand{\Anear}{\Amat^\tbox{near}}
\newcommand{\Afar}{\Amat^\tbox{far}}
\newcommand{\Aproxy}{\Amat^\tbox{proxy}}
\newcommand{\Aself}{\Amat^\tbox{self}}
\newcommand{\Aleft}{\Amat^\tbox{left}}
\newcommand{\Aright}{\Amat^\tbox{right}}

\newcommand{\Tmat}{\mathsf{T}}

\newdefinition{remark}{Remark}[section]
\newdefinition{definition}{Definition}[section]

\journal{J. Comput. Phys.}

\begin{document}
\begin{frontmatter}
\title{A fast direct solver for scattering from periodic structures
       with multiple material interfaces in two dimensions}
\author[lg]{Leslie Greengard}
\address[lg]{Courant Institute,
New York University, New York, NY 10012.}

\author[kh]{Kenneth L. Ho}
\address[kh]{Department of Mathematics, Stanford University, Stanford, CA 94305.}

\author[jylee]{June-Yub Lee}
\address[jylee]{Department of Mathematics,
Ewha Womans University, Seoul, 120-750, Korea.}


\begin{abstract}
We present a new integral equation method for the calculation
of two-dimensional scattering from periodic structures
involving {\em triple-points}
(multiple materials meeting at a single point).
The combination of a robust and high-order accurate
integral representation and a fast direct solver
permits the efficient simulation of scattering from fixed structures
at multiple angles of incidence.
We demonstrate the performance of the scheme with
several numerical examples.
\end{abstract}

\end{frontmatter}

\section{Introduction}

\begin{figure}[!htb]
\begin{center}
\includegraphics[width=4in]{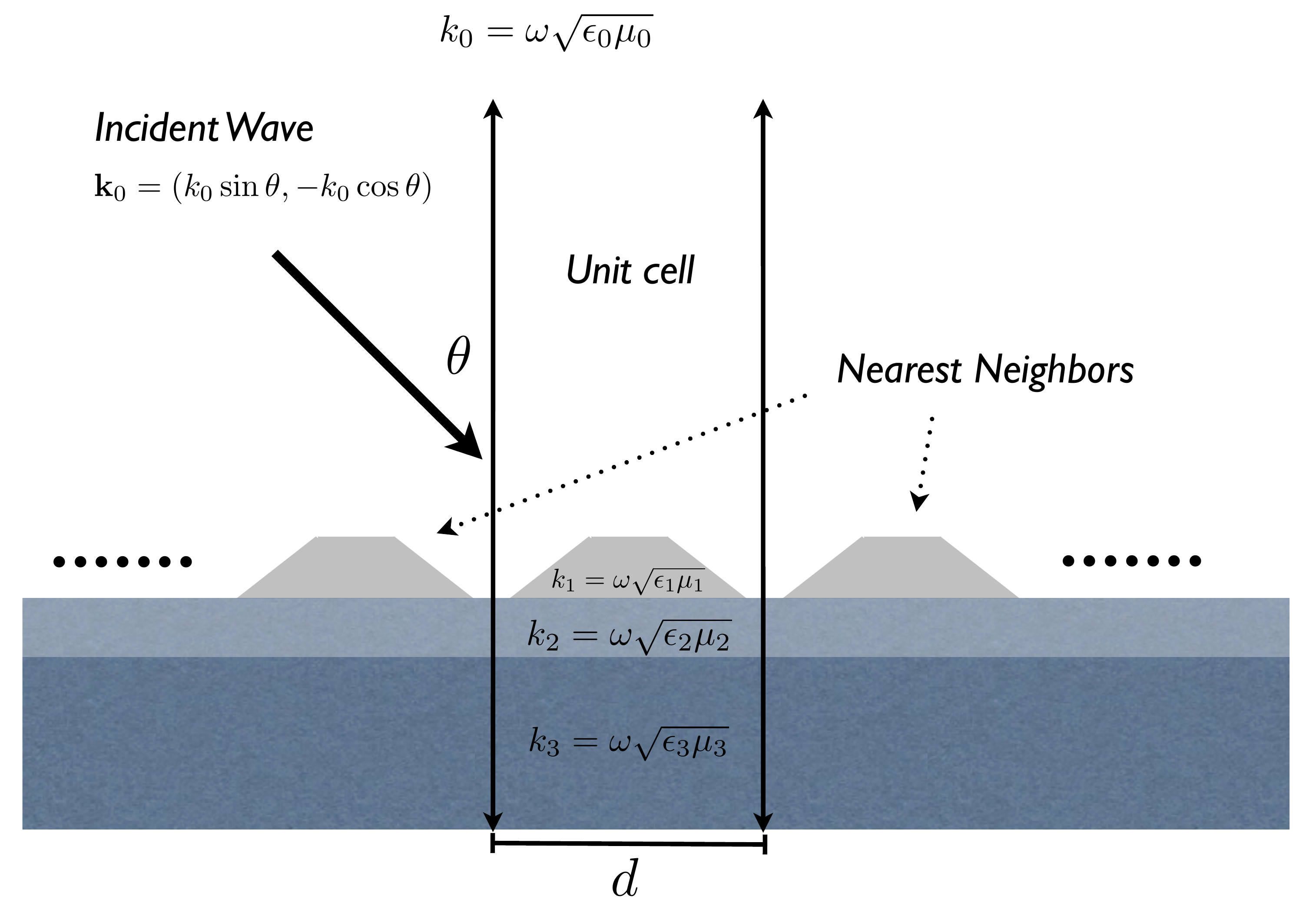}%
\end{center}
\caption{A periodic array of scatterers on the surface of a layered medium.
The Helmholtz coefficient for the upper medium is $k_0$, that for
the trapezoidal-shaped scatterers is $k_1$ and that of the two
layers beneath are $k_2$ and $k_3$, respectively. We assume that
the lowest interface (here between the $k_2$ and $k_3$ layers) is
located at $y=0$ and that the maximum height of the scatterers is at
$y=y_0$. We also assume that the unit cell is centered at $x=0$.
The bottom layer is assumed to be infinite in extent.}
\label{fig-geom1}
\end{figure}

The interaction of acoustic or electromagnetic waves with structured,
periodic materials is often complicated by the fact that the scattering
geometry involves domains where multiple media meet at a single point.
Examples include the design of
diffraction gratings, the development of high efficiency solar cells,
and non-destructive optical inspection in semiconductor manufacturing
(metrology) \cite{ATWATER,CTV,FMMBOOK,METROLOGY,PM,YEUNG,YTS}.
The geometry of a typical scattering problem is shown in
Fig.~\ref{fig-geom1}.

For the sake of concreteness, we will assume throughout this paper that
the governing equations are the Maxwell equations
in two dimensions (here, the $xy$ plane). We also assume the
incident wave is in TM-polarization
\cite{JACKSON} and that each of the constituent materials is
locally isotropic with constant permittivity $\epsilon$ and
permeability $\mu$.
In this case, the Maxwell equations are well-known to take the
simpler form
\begin{align*}
 \E(x,y,z) &= \E(x,y) = \,
(0, 0, E(x,y))\\
 \H(x,y,z) &= \H(x,y) = \, \frac{1}{i \omega \mu}
(E_y(x,y), -E_x(x,y), 0)
\end{align*}
with
\begin{equation}
\label{Helmtot}
\nabla^2 E(\x) + k^2(\x) E(\x) = 0 \quad
{\rm  for }\; \x = (x,y) \in R^2 \, .
\end{equation}
Here,
$k(\x) = \omega \sqrt{\epsilon(\x) \mu(\x)}$, where
we have assumed a time-dependence of $e^{-i\omega t}$ with
$\omega>0$ the frequency of interest.

Using the language of scattering theory, we let
\begin{equation}
\label{phitotdef}
E(\x) = u^\tbox{in}(\x) + u(\x) \, ,
\end{equation}
where $u^\tbox{in}(\x)$ is a known incoming field,
\[ u^\tbox{in}(\x) = u^\tbox{in}_\theta(x,y) = e^{i k_0 (\sin\theta\: x - \cos\theta\: y)}, \]
and $u(\x)$ is the unknown scattered field.
At material interfaces,
\begin{align}
\label{dielecbc1}
\left[E \right] =  0  \quad &\Rightarrow\quad
\left[u \right] = -\left[u^\tbox{in}\right]  \\
\label{dielecbc2}
\left[\frac{1}{\mu} \frac{\partial E}{\partial \nu} \right] =  0
\quad &\Rightarrow\quad
\left[\frac{1}{\mu} \frac{\partial u}{\partial \nu} \right] =
-\left[\frac{1}{\mu} \frac{\partial u^\tbox{in}}{\partial \nu}\right] \, ,
\end{align}
where $\nu$ denotes the normal direction and $\left[ f \right]$ denotes the
jump in the quantity $f$ across an interface.
For simplicity, we will assume $\mu=1$ and $\epsilon$ is distinct in each
domain.
The essential difficulties that we wish to address are manifested
in that setting, so we ignore other variants of the
scattering problem without loss of generality.

Scattering problems of the type illustrated in
Fig.~\ref{fig-geom1} are often called {\em quasi-periodic}
since the obstacles are arrayed periodically, but the incoming,
scattered and total field experience a phase change in traversing the
unit cell:
\begin{equation}
u(x+d,y) = e^{i \alpha d} u(x,y) ,
\label{qpcond}
\end{equation}
where $\alpha = k_0 \sin \theta$. (In this convention, normal incidence
corresponds to $\theta = 0$.)

In the $y$-direction, to obtain a well-posed problem, the scattered field
$u$ must satisfy a somewhat involved
radiation condition \cite{QPSCAT,bonnetBDS,nedelecstarling,shipmanreview}
- namely that it takes the form of Rayleigh--Bloch expansions
\begin{eqnarray}
u(x,y) &=& \sum_{n\in\mathbb{Z}} a^+_n
   e^{i \kappa_n x} e^{i k_n y}
   \qquad y>y_0, \; x\in\mathbb{R}
\label{e:rbu}
\\
u(x,y) &=& \sum_{n\in\mathbb{Z}} a^-_n
   e^{i \kappa_n x} e^{-i k^{(-)}_n y}
   \qquad y<0, \; x\in\mathbb{R}  \, ,
\label{e:rbd}
\end{eqnarray}
assuming, as in Fig.~\ref{fig-geom1}, that the lowest interface
lies at $y=0$ and that $y_0$ is the maximum extent of the scatterers.
In this formula,
$\kappa_n = k_0 \sin \theta + \frac{2\pi n}{d}$,
in order to satisfy the quasi-periodicity condition.
Letting $k_n = +\sqrt{k_0^2-\kappa_n^2}$ enforces that the
expansion satisfy the homogeneous Helmholtz equation in the upper
half-space, while
letting $k^{(-)}_n = +\sqrt{k_l^2-\kappa_n^2}$ enforces that the
expansion satisfy the homogeneous Helmholtz equation in the lower
half-space with wavenumber $k_l$ ($k_3$ in Fig.~\ref{fig-geom1}).

Above the scatterers in the unit cell ($y>y_0$), note that
if $|\kappa_n|\le k_0$, then $k_n$ is real and the
waves in the Rayleigh-Bloch expansion
(\ref{e:rbu}) are {\em propagating} modes.
If $|\kappa_n| > k_0$, then $k_n$ is imaginary and the
corresponding modes are called {\em evanescent}. They
do not contribute to the far field.

\begin{definition}
The complex coefficients $a^+_n$ for propagating modes in the Rayleigh-Bloch
expansion are known as the {\em Bragg diffraction amplitudes at the grating
orders}.
\end{definition}

For each fixed $\alpha$ and $d$, there is a discrete set of
frequencies $\omega$ for which some $k_n$ may vanish, at which
point the Rayleigh-Bloch mode is constant in the $y$-direction.
Such modes are called {\em Wood's anomalies}.
(There is also a discrete set of frequencies where the solution is nonunique,
due to guided modes which propagate along the grating.
The latter are, in a certain sense, nonphysical and
we refer the interested reader to \cite{bonnetBDS,lintonthompson,shipmanreview}
for further discussion.)

In this paper, we present an integral equation method and a corresponding
fast direct solver for scattering problems of the type discussed above.
We make use of the {\em quasi-periodic Green's function}, which requires only
a discretization of the dielectric interfaces within the unit cell.
In a recent paper, Gillman and Barnett
\cite{gillmanbarnett} address the same problem using a slightly different
formulation with a different approach to imposing quasi-periodicity.
We will discuss the relative advantages of the two approaches in
section \ref{sec:conclusions}.

\section{The quasi-periodic Green's function}
\label{s:gqp}

A classical approach to the calculation of
quasi-periodic scattering is based on using the Green's function
that satisfies the desired conditions (\ref{qpcond}), (\ref{e:rbu}), and
(\ref{e:rbd})
\cite{arens,linton,nedelecstarling,nicholas,otani08,venakides00}.
This is accomplished by constructing a
one-dimensional array of suitably ``phased" copies of the free-space
Green's function for the Helmholtz equation with wavenumber $k$.
More precisely, the quasi-periodic Green's function is defined by
\begin{equation}
\gqp(\x) =
\gqp^{(k,\alpha,d)}(\x) = \frac{i}{4}
\sum_{m = -\infty}^\infty e^{im\alpha d} H_0^{(1)}(k|\x-(md,0)|) \, ,
\label{e:gqp}
\end{equation}
where $H_0^{(1)}$ is the outgoing Hankel function of order zero.
It is clear that the sum formally satisfies the condition (\ref{qpcond}).
The Rayleigh-Bloch conditions (\ref{e:rbu}),
(\ref{e:rbd}) follow from Fourier analysis and the fact
that $H_0^{(1)}$ itself satisfies the Sommerfeld radiation condition.
Unfortunately, the series in
(\ref{e:gqp}) is only conditionally convergent for real $k$. To obtain
a physically meaningful limit, one adds a small amount of dissipation
($k \rightarrow k + i \eps$) and considers
$\lim_{\eps \rightarrow 0} \gqp^{(k+i\eps,\alpha,d)}(\x)$.
(See \cite{QPSCAT,dienst01} for a more detailed discussion).
We define the
``near field'' of the quasi-periodic Green's function by
\begin{equation}
\gqpnear(\x) = \frac{i}{4} \sum_{m\in [-1,0,1]}
e^{im\alpha d} H_0^{(1)}(k|\x-(md,0)|)
\label{e:gqpn}
\end{equation}
and the ``smooth'' part of the quasi-periodic Green's function by
\begin{equation}
\gqpfar(\x) =  \frac{i}{4} \sum_{\substack{m\in\mathbb{Z} \\
m \neq [-1,0,1]}} e^{im\alpha d} H_0^{(1)}(k|\x-(md,0)|) \, .
\label{e:gqpr}
\end{equation}
The latter is a smooth solution to the Helmholtz equation
within the unit cell centered at the origin (see Fig.~\ref{fig-geom1})
and can be expanded in a Bessel series
\begin{equation}
\gqpfar(\x) =  \sum_{n=-\infty}^{\infty} s_n J_n(k|\x|).
\label{jexp}
\end{equation}
In the low frequency regime, where the unit cell is on the order of
a few wavelengths or smaller,
the Bessel series converges rapidly so long as
the $y$-component of the target point $\x$ is less than $d$.
For larger values of $y$ it is more convenient to switch representations
and use the Rayleigh-Bloch expansion (\ref{e:rbu}) directly.
An analytic formula for the coefficients
$s_n$ of the Bessel expansion (\ref{jexp}) can be obtained from
the Graf addition theorem \cite[Eq. 9.1.79]{AS}:
\begin{equation}
s_n = \frac{i}{4} \sum_{\substack{m\in\mathbb{Z} \\
m \neq [-1,0,1]}} e^{im\alpha d} H_n(k|md|) (-1)^{n \cdot \rm{signum}(m)}
\end{equation}
These coefficients are known as {\em lattice sums} and
depend only on the parameters $k,\alpha,d$.
Most numerical schemes for the rapid evaluation of the
quasi-periodic Green's function are based on the evaluation of
\begin{equation}
\gqp(\x) = \frac{i}{4} \sum_{m\in [-1,0,1]}
e^{im\alpha d} H_0^{(1)}(k|\x-(md,0)|) +
\sum_{n=-\infty}^{\infty} s_n J_n(k|\x|),
\end{equation}
combining (\ref{e:gqpn}) and (\ref{jexp}).
There is a substantial literature on efficient methods for computing
the lattice sums themselves (see, for example,
\cite{dienst01,linton98,McPh00,moroz01}). In this paper
we use a scheme based on asymptotic analysis and the Euler-MacLaurin
formula \cite{ryan}. Since there are a number of effective schemes for this
step, we omit further discussion except to note that
\begin{enumerate}
\item the
quasi-periodic Green's function fails to exist at Wood's anomalies
\item if the scattering structure in the unit cell has a high aspect
ratio $y_0 \gg d$, then the lattice sum approach is inconvenient because
more images need to be added to $\gqpnear$ in order to ensure convergence
of the Bessel expansion for $\gqpfar$.
\end{enumerate}
We refer the reader to
\cite{QPSCAT,gillmanbarnett} for a method capable of handling both these difficulties.
Here, we assume that $\gqp$ is well-defined and that the aspect ratio
$y_0/d$ is less than or equal to 1.

\section{The integral equation}
\label{s:inteq}

In the absence of triple-points, a number of groups have developed
high-order accurate integral equation methods for scattering from
periodic structures (see, for example,
\cite{arens,QPSCAT,BRUNOHASLAM,HAIDER02,nicholas,YEUNG}).
For this, suppose that we have a single scatterer
$\Omega$ in the unit cell,
with Helmholtz parameter $k_1$ and boundary $\Gamma$.
In the context of Fig.~\ref{fig-geom1}, this would correspond to
an absence of the layered substrate (that is, $k_2=k_3=k_0$), with
$\Omega$ an isolated trapezoidal-shaped scatterer.
One can then use the
representation
\begin{equation}
u \;=\; \left\{\begin{array}{ll}
{\Sc}^{k_1}[\Gamma,\sigma](\x) +
{\Dc}^{k_1}[\Gamma,\mu](\x) & \mbox{in } \Omega\\
{\Sc}^{k_0}_\tbox{QP}[\Gamma,\sigma](\x) +
{\Dc}^{k_0}_\tbox{QP}[\Gamma,\mu](\x)& \mbox{in }
U\setminus\overline{\Omega}
\end{array}\right.
\label{e:qprep}
\end{equation}
where $\Sc$ and $\Dc$ denote the usual single and double layer operators
\cite{CK,nedelec,GLEE}
\begin{eqnarray}
{\Sc}^{k}[\Gamma,\sigma](\x) &=& \int_\Gamma G^{k}(\x-\y) \sigma(\y) ds_\y
\label{e:s}
\\
{\Dc}^{k}[\Gamma,\mu](\x) &=& \int_\Gamma
\frac{\partial G^{k}}{\partial \nu_\y}(\x-\y) \mu(\y) ds_\y \, .
\label{e:d}
\end{eqnarray}
with $G^k(\x) = \frac{i}{4} H_0^{(1)}(k|\x|)$.
The quasi-periodic layer potentials
$\SQP^{k}$ and $\DQP^{k}$ are simply defined
by replacing the
free-space Green's function $G^k(\x)$ with $\gqp^k(\x)$.
Here $ds$ indicates that we are integrating in arclength on $\Gamma$,
and $\nu_{\y}$ denotes the outward normal at $\y \in \Gamma$.
We will also need the normal derivatives of $S^k$ and $D^k$ at a point
$\x \in \Gamma$, defined by
\begin{equation}
 \Nc^{k}(\Gamma,\sigma,\x) = \int_{\Gamma}
\frac{\partial G^{k}}{\partial \nu_{\x}}(\| \x - \y \|) \,\sigma(\y) \, ds_{\y}
\ ,\
 \Tc^{k}(\Gamma,\mu,\x) = \int_{\Gamma}
 \frac{\partial^2 G_{k}}{\partial \nu_{\x}\partial \nu_{\y}}
(\| \x - \y \|) \,\mu(\y) \, ds_{\y}
\, .
\end{equation}
The periodic versions $\NQP$ and $\TQP$ are defined in the same manner.
Note that by construction, the governing Helmholtz equation is satisfied
in each domain.
Note also that we have only chosen to use the quasi-periodic layer potentials
in the exterior domain $U \setminus \Omega$. In the context of Fig.~\ref{fig-geom1}, we will use the quasi-periodic layer potentials for
the $k_0$, $k_2$ and $k_3$ domain and the standard layer potentials
for the $k_1$ domain.
$\Sc^k$ is weakly singular as $\x \rightarrow \Gamma$, and the integral
is well-defined. For $\Dc^{(k)}$ and $\Nc^{(k)}$,
the limiting value depends on the
side of $\Gamma$ from which $\x$ approaches the curve. For $\x \in  \Gamma$,
we assume both are defined in the principal value sense.
The operator $\Tc^{(k)}$ is {\em hypersingular} and unbounded as a map
from the space
of smooth functions on $\Gamma$ to itself. It should be interpreted
in the Hadamard finite part sense.

Substituting the representation
(\ref{e:qprep})
into the interface conditions
(\ref{dielecbc1}), (\ref{dielecbc2})
and taking the appropriate limits yields the system of integral equations
\begin{subequations}\label{grpb}
\begin{align}
\mu(\x)
&+ (\SQP^{k_{0}}(\Gamma,\sigma) - {\Sc}^{k_{1}}(\Gamma,\sigma))[\x] +
(\DQP^{k_0}(\Gamma_,\mu) - \Dc^{k_1}(\Gamma,\mu))[\x] =
- \left[u^\tbox{in}(\x) \right]
\label{IB1}\\
-\sigma(\x) &+
 (\NQP^{k_{0}}(\Gamma,\sigma) - \Nc^{k_{1}}(\Gamma,\sigma))[\x] +
(\TQP^{k_{0}}(\Gamma,\mu) - \Tc^{k_{1}}(\Gamma,\mu))[\x]
=
- \left[ \frac{\partial u^\tbox{in}}{\partial \nu}(\x) \right]
\label{IB2}
\end{align}
\end{subequations}
for the unknowns $[\sigma,\mu]$.

A critical feature of the system (\ref{IB1}), (\ref{IB2}) is that,
while $\Tc$ itself is hypersingular,
only the {\em difference} of hypersingular kernels
appears in the equations. All the operators appearing above are
compact on smooth domains and we have a
system of Fredholm equations of the second kind,
for which the formal theory is classical \cite{GLEE,MIKHLIN} and the
solution is unique.
The cancellation of hypersingular terms in this manner was introduced
in electromagnetics by M\"{u}ller \cite{MULLER}, and in the scalar
case by Kress, Rokhlin, Haider, Shipman and Venakides
\cite{HAIDER02,Kressieq,ROK83}.

For smooth domains, the issue of quadrature
has been satisfactorily resolved, so that
high order accuracy is straightforward
to achieve \cite{ALPERT,BGR,HELSING,HO,QBX,KRESS91}. The generalized Gaussian
quadrature method of \cite{BGR}, for example,
permits the use of composite quadrature rules that take into
account the singularity of the Green's function and can be stored in
tables that do not depend on the curve geometry.
Assuming the boundary component $\Gamma$ is subdivided into $P$
curved panels with
$k$ points on each panel, these rules achieve $k$-th order accuracy.
More precisely, each integral operator
\[ \int_\Gamma G^{k}(\x-\y) \sigma(\y) ds_\y \]
is replaced by a sum of the form
\[ \sum_{p=1}^P \sum_{j=1}^{k} {\cal G}^{k}(x_{q,l},y_{p,j}) \sigma_{p,j}
w_{q,l,p,j}
\]
where $x_{q,l}$ is the $l$-th Gauss-Legendre node on panel $q$,
$y_{p,j}$ is the $j$-th Gauss-Legendre node on panel $p$,
$w_{q,l,p,j}$ is a quadrature weight and
${\cal G}^{k}(x_{q,l},y_{p,j})$ is a ``quadrature kernel".

\begin{remark}
For nonadjacent panels,
${\cal G}^{k}(x_{q,l},y_{p,j})$ is simply the original kernel
$G^{k}(x_{q,l},y_{p,j})$. For the interaction of a panel with itself or
its two nearest neighbors, the quadrature kernel is produced by a somewhat
involved interpolation scheme according to the
generalized Gaussian quadrature formalism \cite{BGR}.
From a linear algebra perspective, generalized Gaussian
quadrature can be viewed as producing a block tridiagonal matrix
(with block size $k \times k$)
of interactions of each panel with
itself and its two neighbors. These are computed directly.
All other block matrix interactions are obtained using
standard Gauss-Legendre weights $w_{q,l,p,j} = w_{p,j}$ scaled
to the dimensions of the $p$-th source panel.
This structure of the far-field interactions permits the straightforward
use of fast multipole acceleration and the hierarchical direct solver
of \cite{HG}.
\end{remark}

In domains with corners, but not multi--material junctions,
exponentially adaptive grids maintain high-order
accuracy (see, for example \cite{BRS,HO}).
In the simplest version, one can first divide the boundary into
equal size subintervals and employ a $k$-th order generalized Gaussian quadrature rule on each.
For each segment that impinges on a corner point,
one can further subdivide it using a
dyadically refined mesh, creating $\log_2 (1 / \varepsilon)$ additional
subintervals, where $\varepsilon > 0$ is a specified numerical precision. If the same
$k$-th order rule is used for each refined subinterval, it is
straightforward to show that the resulting rule
has a net error of the order
$O\left( e^{-k} \log_2 (1 / \varepsilon) \right)$.
The need for dyadic refinement
comes from the fact that the densities $\sigma$ or $\mu$
may develop singularities at the corner points and the refinement yields a
high order piecewise polynomial approximation of the density.
For $\varepsilon = 10^{-14}$ and $k = 16$, the net corner error
is about $10^{-14}$ while for $k=8$, it is about $10^{-8}$
(see Fig.~\ref{figskel} for an illustration).

\begin{remark}
In recent work, \cite{BRS,HO} have shown that
one can dramatically reduce the
number of degrees of freedom in the vicinity of the corner by the
use of compression, as well. We have not used such optimization here.
\end{remark}

It is now appreciated (see, for example, \cite{BREMER,HELSING})
that the condition number of a properly discretized system of equations
is very well controlled.
Following discretization, we
use Bremer's approach \cite{BREMER} here, which involves setting
the discrete variables to be $\sigma_i \sqrt{w_i}$ and
$\mu_i \sqrt{w_i}$, rather than the density values $\sigma_i$ and $\mu_i$
themselves. This ensures that that spectrum of the discrete system
approximates the spectrum of the continuous integral equation in $L^2$.
The formal analysis is somewhat involved, since operators that are
compact on smooth domains are only bounded (but not compact) on domains with
corners. We refer the reader to
\cite{BREMER,HELSING} for details.

\section{Stable and accurate integral formulations in the presence of
multi--material junctions}

In the case of multiple subdomains,
a natural approach would be to represent the
field in each subdomain $\Omega_i$ with Helmholtz coefficient
$k_i$ in terms of layer potentials on the boundary $\Gamma_i$ of
$\Omega_i$. That is, in subdomain $\Omega_i$, we would represent the
solution as
\begin{equation}
u_i(\x) \;=\;
{\Sc}^{k_i}[\Gamma_i,\sigma](\x) +
{\Dc}^{k_i}[\Gamma_i,\mu](\x)
\label{e:qprepi}
\end{equation}
with $\Sc$ and $\Dc$ replaced by their quasi-periodic counterparts
for subdomains that extend across the unit cell
(the $k_0$, $k_2$, and $k_3$  domains in Fig.~\ref{fig-geom1}).

In doing so, it turns out that
the analog of equations (\ref{IB1},\ref{IB2}) fails to converge
in the presence of multi-material junctions.
The reason for this is simple, and analyzed in \cite{triplepoint}.
Consider the interface condition (\ref{IB2}) for $\xy$ lying
on the segment $\overline{AB}$ in Fig.~\ref{figskel}.
Restricting our attention just to the
segments impinging on the corner point $A$,
we have
\begin{align}
-\sigma(\x) +
\left[ \NQP^{k_{0}}(\overline{AB},\sigma,\xy) -
 \Nc^{k_{1}}(\overline{AB},\sigma,\xy)\right] +
 \NQP^{k_{0}}(\overline{AE},\sigma,\xy) -
 \Nc^{k_{1}}(\overline{AD},\sigma,\xy) \ +\ \dots \hspace{.3in}
 \nonumber \\ +
\left[ \TQP^{k_{0}}(\overline{AB},\mu,\xy) -
 \Tc^{k_{1}}(\overline{AB},\mu,\xy) \right] +
 \TQP^{k_{0}}(\overline{AE},\mu,\xy) -
 \Tc^{k_{1}}(\overline{AD},\mu,\xy) \ +\ \dots
=
- \left[ \frac{\partial u^\tbox{in}}{\partial \nu}(\x) \right]
\label{IB2a}
\end{align}
Note that both the terms $\TQP^{k_{0}}(\overline{AE},\mu,\xy)$ and
$\Tc^{k_{1}}(\overline{AD},\mu,\xy)$ involve hypersingular contributions
at the junction $A$ without forming part of a difference kernel.
This destroys the high-order accuracy of the scheme.

By using a global integral representation, it was shown in
\cite{triplepoint} that high-order accuracy can be restored.
That is, instead of (\ref{e:qprepi}), we let
\begin{equation}
u_i(\x) \;=\;
{\Sc}^{k_i}[\Gamma,\sigma](\x) +
{\Dc}^{k_i}[\Gamma,\mu](\x)
\label{e:qprepall}
\end{equation}
and apply the continuity conditions.
For $\xy \in \Gamma$ lying on an interface between subdomains
with Helmholtz coefficients $k_i$ and $k_j$, we have
\begin{subequations}\label{grpc}
\begin{align}
\mu(\xy) &+
 \Sc^{k_{i}}(\Gamma,\sigma,\xy)  -
 \Sc^{k_{j}}(\Gamma,\sigma,\xy) + 
 \Dc^{k_{i}}(\Gamma,\mu,\xy)  -
 \Dc^{k_{j}}(\Gamma,\mu,\xy) =
- \left[\Phi^\tbox{in}(\xy) \right]
\label{IC1}\\
- \sigma(\xy) &+
 \Nc^{k_{i}}(\Gamma,\sigma,\xy) -
 \Nc^{k_{j}}(\Gamma,\sigma,\xy) + 
 \Tc^{k_{i}}(\Gamma,\mu,\xy) -
 \Tc^{k_{j}}(\Gamma,\mu,\xy)
=
- \left[ \frac{\partial \Phi^\tbox{in}}{\partial \nu}(\xy) \right]
\label{IC2}
\end{align}
\end{subequations}
As above, the operators $\Sc$, $\Dc$, $\Nc$, $\Tc$ are
replaced by their quasi-periodic counterparts
for subdomains that extend across the unit cell
(the $k_0$, $k_2$, and $k_3$  domains in Fig.~\ref{fig-geom1}).

The global representation (\ref{e:qprepall}) is ``non-physical"
in the sense that the field in a given subdomain is determined, in part,
by layer potential components that are not actually
part of the subdomain's boundary. By doing so, however,
we remove all hypersingular terms from the integral equation.
Only difference kernels appear in the final linear system.
One could improve efficiency somewhat, while achieving similar results,
by supplementing the representation (\ref{e:qprepi})
only by the boundary segments that actually impinge on a multi-material
junction. We use the fully global representation in our experiments here
for the sake of simplicity.

\begin{remark}
For related approaches addressed at solving problems with multi-material
junctions, see
\cite{CLAEYS,CLAEYS2,HJH}.
\end{remark}

\section{Fast direct solvers}

Given a well-conditioned and high order discretization,
large scale scattering problems in singular geometries
can be solved by using fast multipole-accelerated
iterative solution methods such as GMRES \cite{GMRES}.
While these are asymptotically optimal schemes, one is often
interested in modeling
the interaction of a given physical structure
(such as the geometry in Fig.~\ref{fig-geom1})
with a large number of incoming fields. This requires the solution
of an integral equation with multiple right-hand sides, and
standard iterative methods do not take maximal advantage of this fact.

Direct solvers, on the other hand, first construct a factorization
of the system matrix, then solve against each right-hand side using that factorization at a cost that is typically
much lower. In the last decade, specialized versions have
been created which are particularly suited to the integral equation
environment. This is an active area of research and we do not seek to review
the literature, except to note selected important
developments in the case of hierarchically semiseparable matrices
\cite{chandrasekaran:2006b:siam-j-matrix-anal-appl,chandrasekaran:2006a:siam-j-matrix-anal-appl,xia:2009:siam-j-matrix-anal-appl},
$\mathcal{H}$-matrices
\cite{hackbusch:1999:computing,hackbusch:2002:computing,hackbusch:2000:computing},
and hierarchically block separable matrices
\cite{gillman:2012,greengard:2009:acta-numer,HG,martinsson:2005:j-comput-phys}.
We provide a brief description of the approach,
following the presentation of \cite{gillman:2012,HG}.

\subsection{Recursive skeletonization for integral equations}
Let $\Amat \in \mathbb{C}^{N \times N}$ be the matrix discretization of an integral equation such as (\ref{grpc}), and let its indices $\{ 1, \dots, N \}$ be ordered hierarchically according to a quadtree on the unit cell. This can be done by first enclosing the set of all associated points within a sufficiently large box. If the box contains more than a specified number of points, it is subdivided into four quadrants and its points distributed accordingly between them. This procedure is repeated for each new box added, terminating only when all boxes contain $O(1)$ points. The boxes that are not subdivided are called {\em leaf boxes}. For simplicity, we assume that all leaf boxes live on the same level of the tree, but this restriction can easily be relaxed.

Start at the bottom of the tree and consider the partitioning induced by the leaves. Let $p$ be the number of leaf boxes and assume that each has $n$ points so that $N = pn$. Then $\Amat$ has the block form $\Amat = \Amat_{ij}$ for $i, j = 1, \dots, p$. We now use the {\em interpolative decomposition} (ID) \cite{cheng:2005:siam-j-sci-comput} to skeletonize $\Amat$. The ID is a matrix factorization that rewrites a given low-rank matrix in terms of a subset of its rows or columns, called {\em skeletons}. In the integral equation setting, the off-diagonal block rows
\begin{equation}
 \Amat_{i \leftarrow} = [\Amat_{i,1}, \: \dots, \: \Amat_{i,i-1}, \: \Amat_{i,i+1}, \: \dots, \: \Amat_{i,p}]
 \label{eq:block-row}
\end{equation}
are low-rank due to the smoothness of the Green's function, and the same is true of the off-diagonal block columns. Thus, it can be shown \cite{gillman:2012,HG} that the ID enables a representation of the form
\begin{equation}
 \Amat_{ij} = \Lmat_{i} \Smat_{ij} \Rmat_{j}, \quad i \neq j,
 \label{eq:block-separable}
\end{equation}
for each off-diagonal block, where $\Lmat_{i} \in \mathbb{C}^{n \times k}$, $\Rmat_{j} \in \mathbb{C}^{k \times n}$, and $\Smat_{ij} \in \mathbb{C}^{k \times k}$ is a submatrix of $\Amat_{ij}$, with $k \ll n$. The matrix can then be written as
\begin{equation}
 \Amat = \Dmat + \Lmat \Smat \Rmat,
 \label{eq:compressed-representation}
\end{equation}
where
$$
\Dmat = \left[
\begin{array}{ccc}
\Amat_{11} & & 0 \\
& \ddots\\
0 & & \Amat_{pp}
\end{array} \right] \in \mathbb{C}^{N \times N}  \, ,
$$
$$
\Lmat = \left[
\begin{array}{ccc}
\Lmat_{1} & & 0 \\
& \ddots\\
0 & & \Lmat_{p}
\end{array} \right] \in \mathbb{C}^{N \times K}, \qquad \Rmat = \left[
\begin{array}{ccc}
\Rmat_{1} & & 0 \\
& \ddots\\
0 & & \Rmat_{p}
\end{array} \right] \in \mathbb{C}^{K \times N}
$$
are block diagonal with $K = pk$, and
$$
\Smat = \left[
\begin{array}{cccc}
0 & \Smat_{12} & \cdots & \Smat_{1p}\\
\Smat_{21} & 0 & \cdots & \Smat_{2p}\\
\vdots & \vdots & \ddots & \vdots\\
\Smat_{p1} & \Smat_{p2} & \cdots & 0\\
\end{array} \right] \in \mathbb{C}^{K \times K}
$$
is dense with zero diagonal blocks.

\begin{remark}
 \label{rmk:proxy}
 The efficient calculation of the interpolation matrices $\Lmat_{i}$ and $\Rmat_{j}$, and the associated skeleton indices, in (\ref{eq:block-separable}) is somewhat subtle. Briefly, it involves separating out neighboring and far-field interactions and representing the latter via free-space interactions with a local ``proxy'' surface. This is justified by the observation that any well-separated interaction governed by a homogeneous partial differential equation (here, the Helmholtz equation) can be induced by sources/targets on the proxy surface, each of which is expressed
in terms of the free-space kernel. For details, see \cite{gillman:2012,HG}. In this paper, for a box of scaled size $1$, we use the circle of radius $1$ about the box center as its proxy surface. Note that all neighbors are defined relative to the periodicity of the unit cell.
\end{remark}

Now consider the linear system $\Amat \xmat = \bmat$. One way to solve it is to construct $\Amat^{-1}$ directly from (\ref{eq:compressed-representation}) using a variant of the Sherman-Morrison-Woodbury formula. This approach is taken in \cite{gillman:2012,martinsson:2005:j-comput-phys}. Here, we follow the strategy of \cite{HG} instead and let $\zmat = \Rmat \xmat$ and $\ymat = \Smat \zmat$ to obtain the equivalent {\em sparse} system
\begin{equation}
 \left[
 \begin{array}{ccc}
  \Dmat & \Lmat\\
  \Rmat & & -\Imat\\
  & -\Imat & \Smat
 \end{array} \right] \left[
 \begin{array}{c}
  \xmat\\
  \ymat\\
  \zmat
 \end{array} \right] = \left[
 \begin{array}{c}
  \bmat\\
  0\\
  0
 \end{array} \right].
 \label{eq:sparse-embedding}
\end{equation}
This can be solved efficiently using any standard sparse direct solver and may provide better stability. In this paper, we use the open-source software package UMFPACK \cite{davis:2004:acm-trans-math-softw,davis:1997:siam-j-matrix-anal-appl}.

Since $\Smat$ is a submatrix of $\Amat$ (up to diagonal modifications), $\Smat$ can itself be expressed in the form (\ref{eq:compressed-representation}) by moving up one level in the tree and regrouping appropriately. This leads to a multilevel decomposition
\begin{equation}
 \Amat = \Dmat^{(\lambda)} + \Lmat^{(\lambda)} \left( \cdots \Dmat^{(1)} + \Lmat^{(1)} \Dmat^{(0)} \Rmat^{(1)} \cdots \right) \Rmat^{(\lambda)},
 \label{eq:multilevel-representation}
\end{equation}
where the superscript indexes the tree level $l = 0, 1, \dots, \lambda$ with $l = 0$ denoting the root. We call this process {\em recursive skeletonization}. The analogue of (\ref{eq:sparse-embedding}) is
\begin{equation}
 \left[
 \begin{array}{cccccc}
  \Dmat^{(\lambda)} & \Lmat^{(\lambda)}\\
  \Rmat^{(\lambda)} & & -\Imat\\
  & -\Imat & \ddots & \ddots\\
  & & \ddots & \Dmat^{(1)} & \Lmat^{(1)}\\
  & & & \Rmat^{(1)} & & -\Imat\\
  & & & & -\Imat & \Dmat^{(0)}
 \end{array} \right] \left[
 \begin{array}{c}
  \xmat^{(\lambda)}\\
  \ymat^{(\lambda)}\\
  \vdots\\
  \xmat^{(1)}\\
  \ymat^{(1)}\\
  \xmat^{(0)}
 \end{array} \right] = \left[
 \begin{array}{c}
  \bmat\\
  0\\
  \vdots\\
  0\\
  0\\
  0
 \end{array} \right],
 \label{eq:multilevel-embedding}
\end{equation}
corresponding to expanding $\Smat$ out in the same way. It can be shown that the solver requires $O(N \log N)$ work when the unit cell is a moderate number of wavelengths in size. We refer the reader to
\cite{gillman:2012,HG} for further discussion.

For our present purposes, we simply note that the output of the fast direct solver is a {\em compressed} representation of the inverse which is computed in two steps:
\begin{enumerate}
 \item
  a recursive skeletonization procedure to obtain the compressed forward operator (\ref{eq:multilevel-representation}); and
 \item
  a factorization of the sparse matrix embedding in (\ref{eq:multilevel-embedding}).
\end{enumerate}
Both steps have the same asymptotic complexity, but the constant for compression is typically far larger. After the inverse has been computed, it can be applied to each right-hand side as needed at a much lower cost.

\begin{remark}
 The ID can be constructed to any specified relative precision $\varepsilon > 0$. This is an input parameter to recursive skeletonization and hence to the direct solver. It can be shown that if (\ref{eq:multilevel-representation}) has relative error $O(\varepsilon)$, as is often the case numerically, then the algorithm produces a solution with relative error $O(\kappa (\Amat) \varepsilon)$, where $\kappa (\Amat)$ is the condition number of $\Amat$. In particular, if $\kappa (\Amat) = O(1)$, as for the integral equation (\ref{grpc}), then the error is $O(\varepsilon)$.
\end{remark}

\begin{remark}
 \label{rmk:in-out}
 Although we have assumed in the discussion above that each block at the same level has the same size, this is in no way essential to the algorithm. In fact, our code uses separate ``incoming'' (row) and ``outgoing'' (column) skeletons for each box. This enables some additional optimization, which, for the present case, can be especially pronounced. This is because while each point receives incoming interactions from only the two wavenumbers on either side of the segment to which it belongs, it sends outgoing interactions consisting of all wavenumbers in the problem. For example, for a point on the segment $\overline{AB}$ in Fig.~\ref{figskel}, it receives at wavenumbers $k_{0}$ and $k_{1}$ but sends at wavenumbers $k_{0}$, $k_{1}$, and $k_{2}$. Therefore, the outgoing skeleton dimension is typically larger, and the amount by which it is larger increases with the total number of wavenumbers/domains.
\end{remark}

\subsection{Multiple angles of incidence}  \label{multangle}
The fast direct solver of the previous subsection allows the robust and accurate solution of
\begin{equation}
 \Amat (\theta) \xmat (\theta) = \bmat (\theta),
 \label{eq:system-angle}
\end{equation}
where we have made explicit the dependence of the integral equation (\ref{grpc}) on the incident angle $\theta$. In the present setting, we are interested in solving (\ref{eq:system-angle}) for many $\theta$. This is not a situation that the direct solver can easily handle since $\Amat (\theta)$ is not fixed. In this subsection, we describe a modified strategy for computing a compressed representation (\ref{eq:multilevel-representation}) of $\Amat (\theta)$ such that it can be rapidly updated to yield a compressed representation of $\Amat (\theta')$ for any $\theta'$ {\em without having to re-skeletonize}. Since skeletonization is typically the most expensive step, this can offer significant computational savings. The sparse matrix in (\ref{eq:multilevel-embedding}) must still be updated and re-factored, but the relative cost of this is small.

To see why such a uniform skeletonization might be possible, consider any finite truncation of the periodic geometry so that it consists merely of a very large array of many, many scatterers. Then the governing integral equation is specified in terms of the {\em free-space} Green's function so that $\Amat$ is independent of $\theta$. The only angle dependence comes from the incoming data $\bmat (\theta)$. Therefore, only one skeletonized representation of $\Amat$ is needed for all $\theta$. The same is true of any finite approximation to the periodic problem.

We now make this intuition precise by considering all interactions, say, incoming on a given box. This is given by the off-diagonal block row (\ref{eq:block-row}) and can be decomposed as
\[
 \Amat_{i \leftarrow} (\theta) = \Anear_{i \leftarrow} (\theta) + \Afar_{i \leftarrow} (\theta)
\]
in terms of the near- and far-field contributions, respectively, to the quasi-periodic Green's function
\[
 \gqp (\x; \theta) = \gqpnear (\x; \theta) + \gqpfar (\x; \theta),
\]
following section \ref{s:gqp}. Clearly, an interpolation basis for both terms together provides an interpolation basis for the sum, so $\Amat_{i \leftarrow} (\theta)$ can be skeletonized by applying the ID to the rows of the matrix
\[
 \tilde{\Amat}_{i \leftarrow} (\theta) = [\Anear_{i \leftarrow} (\theta), \: \Afar_{i \leftarrow} (\theta)].
\]
Since $\gqpfar$ consists only of well-separated interactions, by Remark \ref{rmk:proxy}, $\Afar_{i \leftarrow} (\theta)$ can be replaced by a matrix $\Aproxy_{i \leftarrow}$ corresponding to free-space interactions with a proxy surface. In linear algebraic terms, this means that $\Afar_{i \leftarrow}$ can be written as $\Afar_{i \leftarrow} = \Aproxy_{i \leftarrow} \Tmat_{i \leftarrow} (\theta)$ for some matrix $\Tmat_{i \leftarrow} (\theta)$. Hence,
\begin{equation}
 \tilde{\Amat}_{i \leftarrow} (\theta) = \left[
 \begin{array}{cc}
  \Anear_{i \leftarrow} (\theta) & \Aproxy_{i \leftarrow}
 \end{array} \right] \left[
 \begin{array}{cc}
  \Imat\\
  & \Tmat_{i \leftarrow} (\theta)
 \end{array} \right],
 \label{eq:a-proxy}
\end{equation}
so $\tilde{\Amat}_{i \leftarrow} (\theta)$ can be skeletonized by applying the ID to just the left matrix on the right-hand side. Observe that the angular dependence of the far field has been eliminated.

To eliminate the angular dependence of the near field, we can similarly expand $\gqpnear$ in terms of a $\theta$-independent basis. This can be done using the functions
\[
 \frac{i}{4} H_{0}^{(1)} (k|\x|), \quad \frac{i}{4} H_{0}^{(1)} (k|\x-(d,0)|), \quad \frac{i}{4} H_{0}^{(1)} (k|\x+(d,0)|)
\]
corresponding to interactions with the self-, left-, and right-images, respectively, with corresponding matrices $\Aself$, $\Aleft_{i \leftarrow}$, and $\Aright_{i \leftarrow}$. Then, from (\ref{e:gqpn}),
\begin{equation}
 \Anear_{i \leftarrow} (\theta) = \left[
 \begin{array}{ccc}
  \Aself & \Aleft_{i \leftarrow} & \Aright_{i \leftarrow}
 \end{array} \right] \left[
 \begin{array}{c}
  \Imat\\
  e^{-i \alpha d} \Imat\\
  e^{i \alpha d} \Imat
 \end{array} \right],
 \label{eq:a-near}
\end{equation}
where, recall, $\alpha = k \sin \theta$, so $\Amat_{i \leftarrow} (\theta)$ can be skeletonized by applying the ID to
\begin{equation}
 \tilde{\Amat}_{i \leftarrow} = [\Aself, \: \Aleft_{i \leftarrow}, \: \Aright_{i \leftarrow}, \: \Aproxy_{i \leftarrow}],
 \label{eq:id-matrix}
\end{equation}
which we note has no angular dependence. Thus, the interpolation matrices and skeleton indices resulting from compressing (\ref{eq:id-matrix}) are valid for all $\theta$.

The same approach can be used for outgoing interactions and for interactions at each wavenumber. The result is a modified compressed representation
\begin{equation}
 \Amat (\theta) = \Dmat^{(\lambda)} (\theta) + \Lmat^{(\lambda)} \left( \cdots \Dmat^{(1)} (\theta) + \Lmat^{(1)} \Dmat^{(0)} (\theta) \Rmat^{(1)} \cdots \right) \Rmat^{(\lambda)},
 \label{eq:angle-representation}
\end{equation}
where only the $\Dmat^{(l)} (\theta)$ depend on $\theta$. Therefore, to obtain a compressed representation of $\Amat (\theta')$ for any other $\theta'$, it suffices to perform the update $\Dmat^{(l)} (\theta) \mapsto \Dmat^{(l)} (\theta')$ for each $l$. This, in general, consists only of generating a very small subset of entries of the new matrix and requires $O(N \log N)$ work with a
small constant.

In summary, the full algorithm for analyzing multiple incident angles with fast updating is:
\begin{enumerate}
 \item
  Compress the matrix $\Amat (\theta)$ for some initial $\theta$ by representing interactions with an angle-independent basis such as (\ref{eq:id-matrix}). This is an expensive recursive skeletonization that only needs to be performed once.
 \item
  Embed the resulting decomposition (\ref{eq:angle-representation}) into the sparse matrix of (\ref{eq:multilevel-embedding}) and solve.
 \item
  For each new angle $\theta'$, update the compressed representation (\ref{eq:angle-representation}) via $D^{(l)} (\theta) \mapsto D^{(l)} (\theta')$. Repeat step 2.
\end{enumerate}

\begin{remark}
 In our tests, we have often found it unnecessary to decompose $\Anear_{i \leftarrow} (\theta)$ as in (\ref{eq:a-near}). Instead, we apply the ID to the left matrix on the right-hand side of (\ref{eq:a-proxy}), which depends on $\theta$ but seems to yield results that recover angle independence numerically. This optimization can reduce the constant associated with skeletonization by about a factor of $2$.
\end{remark}

\section{Numerical results} \label{examples}

The algorithm presented above has been implemented in Fortran.
Each boundary segment (in the piecewise smooth boundary) is first divided
into $22$ equal subintervals.
The first and last intervals are then further subdivided with
dyadic refinement toward the corner using $20$ subintervals each.
Thus, the total number of intervals on each
smooth component of the boundary (each side) is $60$ and the number of points
is $480$. We use the 8th order generalized Gaussian quadrature rule of
\cite{BGR} for logarithmic singularities and
solve the integral equations (\ref{grpc})
using recursive skeletonization \cite{gillman:2012,HG} with a tolerance of $\varepsilon = 10^{-9}$.
All timing listed below are for a laptop with a 1.7GHz Intel Core i5
processor.

\begin{figure}[!htb]   
\begin{center}
\includegraphics[width=4.5in]{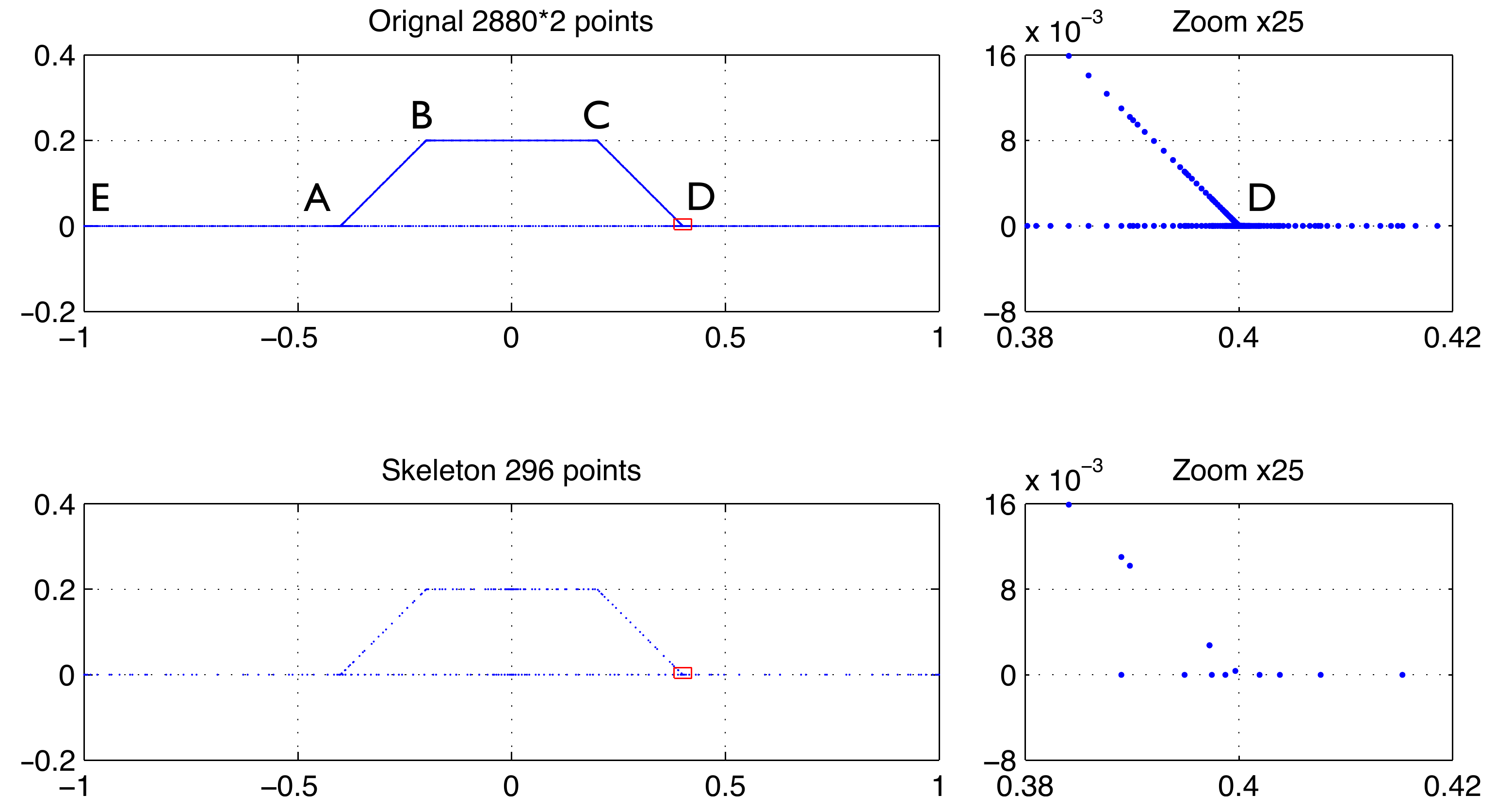}
\end{center}
\caption{Discretization points and skeletons for a quasi-periodic scattering
problem, where the unit
cell consists of a trapezoidal scatterer lying on an infinite substrate.
There are six segments with $480$ points each and $2$ degrees of freedom per
point, resulting in a complex matrix of size $5760 \times 5760$.
After skeletonization, there are only $296$ points remaining at the
coarsest level of the recursion.
\label{figskel}}
\end{figure}

\vspace{.1in}

\begin{figure}[!htb]   
\begin{center}
\includegraphics[width=4.5in]{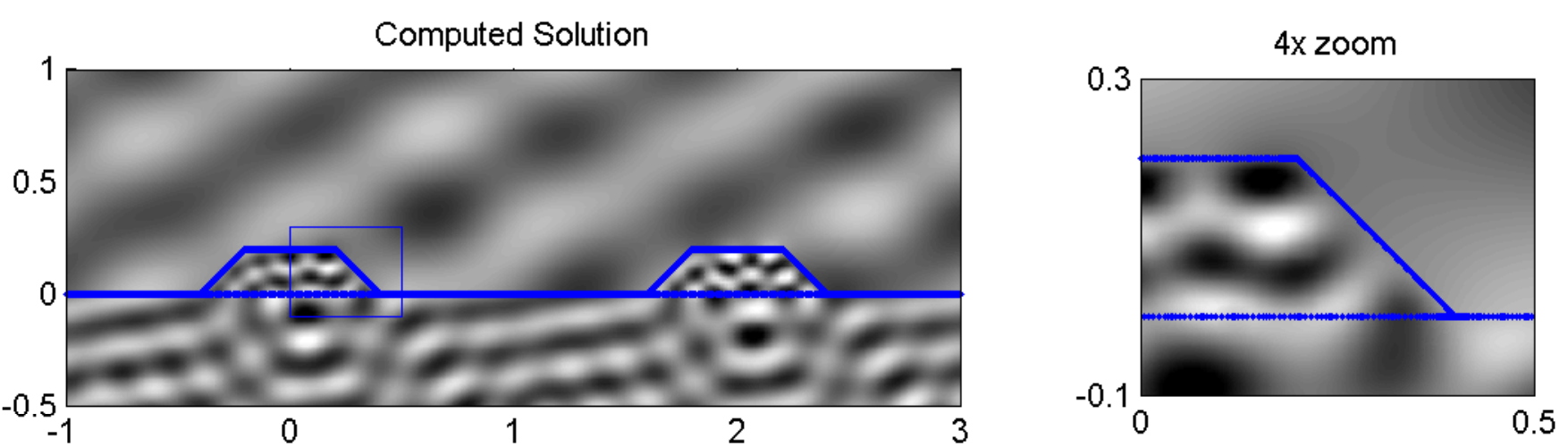}
\end{center}
\caption{The real part of the total field when a plane
wave at $30^{\circ}$ incidence impinges on a periodic structure.
The unit cell and its first neighbor are shown, with a zoom in the
region of the triple-point.}
\label{figs1}
\end{figure}

\noindent
{\bf Example 1:}
We set $\omega=10$, with $\epsilon$ chosen so that the
Helmholtz coefficient in the upper half-space, the trapezoidal
scatterer, and the substrate are $k = 10$, $40\sqrt{2}$, and $30$,
respectively. The incident angle is $30^{\circ}$. The original matrix of dimension $5760 \times 5760$ is compressed to
one of dimension $296 \times 309$. The incoming
and outgoing skeleton dimensions are slightly different as explained in Remark \ref{rmk:in-out} and computed
as part of the recursion.
The time for compression in our current implementation was
290 secs.\ (while generating the necessary matrix entries required 1219 secs.).
Given the compressed representation, the solution time was 2.46 secs.
The resulting accuracy was approximately $10^{-9}$. We plot the
real part of the total field in Fig.~\ref{figs1}.
In Fig.~\ref{figskel}, we plot both the original set of discretization points
and the skeletons that remain at the coarsest level of the recursion.

\medskip \noindent
\subsection{Computing the outgoing modes}

Given our integral representation of the scattered field,
it is straightforward to compute the coefficients
$a^+_n$ in (\ref{e:rbu}) or (\ref{e:rbd}) -
the Bragg diffraction amplitudes at the grating
orders. For an incident field
\[ u^\tbox{in}(\x) = u^\tbox{in}_\theta(x,y) = e^{i k_0 (\sin\theta\: x - \cos\theta\: y)}, \]
we simply let $y_0+\delta$ denote some height
above the scatterers and rewrite
(\ref{e:rbu}) in the form
\[ u(x,y_0+ \delta) e^{-k_0 \sin \theta x} = \sum_{n\in\mathbb{Z}} a^+_n
   e^{2 \pi i nx/d} e^{i k_n (y_0+\delta)} \, ,
\]
where $k_n = +\sqrt{k_0^2-(k_0 \sin \theta + 2\pi n/d)^2}$.
Thus, the $\{ a^+_n \}$ can be computed using Fourier analysis:
\[ a^+_n = \frac{1}{d \, e^{i k_n (y_0+\delta)}}  \int_0^d
u(x,y_0+ \delta) e^{-k_0 \sin \theta x} e^{-2 \pi i n x/d} \, dx \, .
\]
The accurate calculation of $a^+_n$ from this formula
depends on ensuring that the discretization in $x$ is sufficiently
fine to resolve the integrand. In the near field (when $\delta$ is small),
the evanescent modes, corresponding to large $n$, are still present in
$u(x,y_0+ \delta)$ requiring a large number of points to avoid aliasing
errors. By making $\delta$ sufficiently large, the evanescent modes are suppressed.
and a mesh can be used that resolves only the propagating modes - that is,
values of $n$ for which $(k_0 \sin \theta + 2\pi n/d)^2 < k_0^2$.

\vspace{.1in}

\begin{figure}[!htb]   
\begin{center}
\includegraphics[width=4.5in]{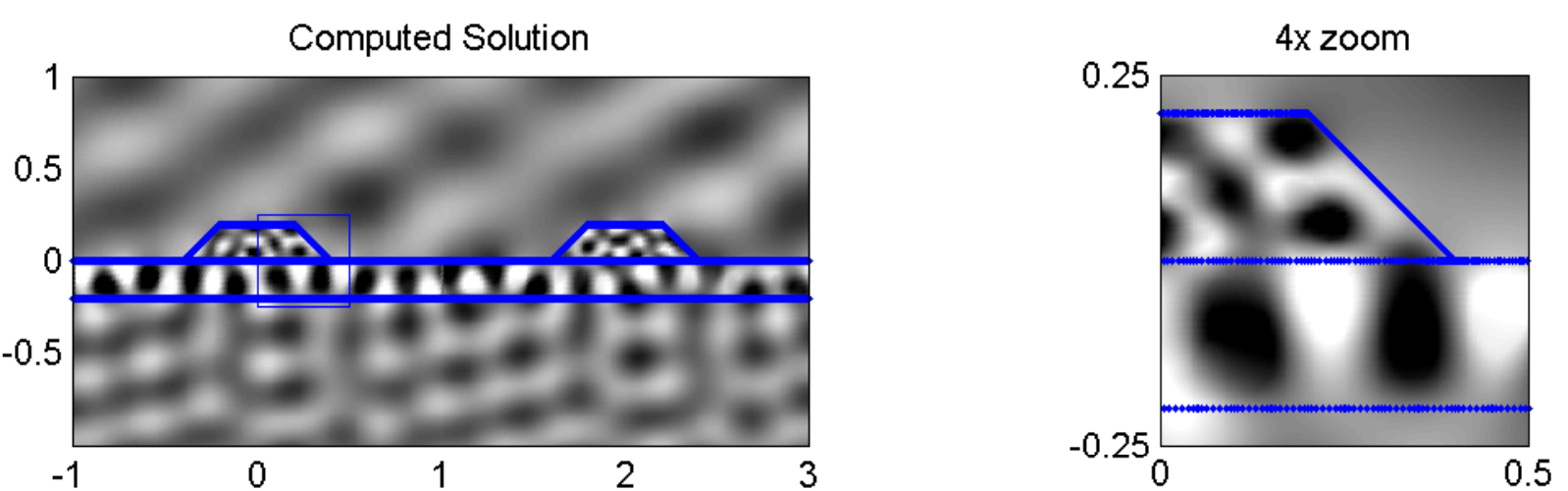}
\end{center}
\caption{ The real part of the total field when a plane
wave at $45^{\circ}$ incidence impinges on a periodic structure with a two-layer
substrate. The unit cell and its first neighbor are shown, with a zoom in the
region of the triple-point.
\label{fig2layer}}
\end{figure}

\begin{figure}[!htb]   
\begin{center}
\includegraphics[width=3.5in]{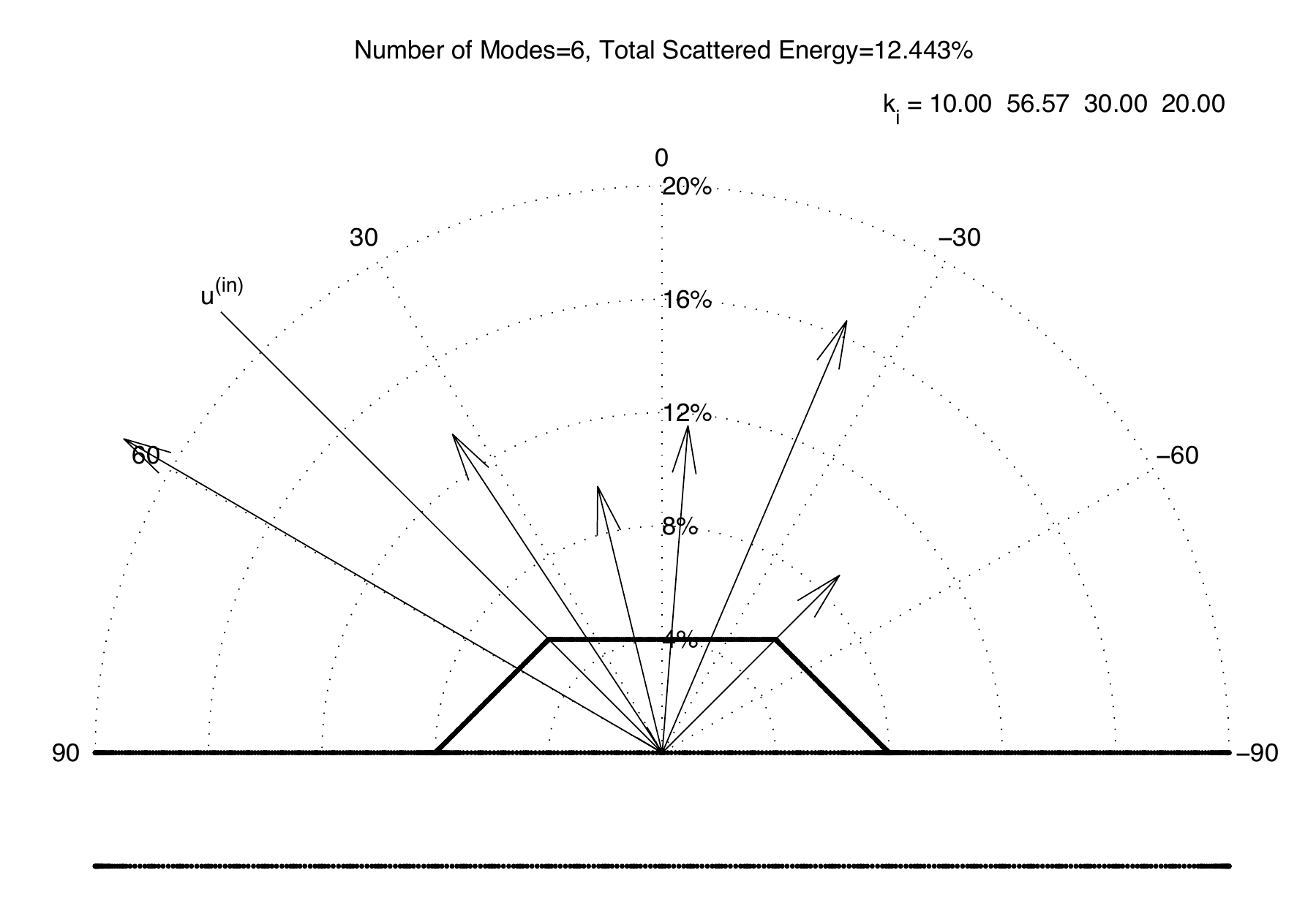}
\end{center}
\caption{The strength of the 6 radiation modes in Example 2 corresponding
to Fig.~\ref{fig2layer}.
The length of the arrows in the various diffraction directions indicate the relative
magnitude of the amplitudes $a^+_n$. Approximately
$12.443 \%$ of the energy is scattered upward.
\label{scatmodes2}}
\end{figure}

\noindent
{\bf Example 2:}
We now consider a scattering problem with
a two-layered substrate (Fig.~\ref{fig2layer}).
We again set $\omega=10$ and choose $\epsilon$ so that the
Helmholtz coefficient in the upper half-space, the trapezoidal
scatterer, and the two substrate layers are $k = 10$, $40\sqrt{2}$, $30$ and $20$,
respectively.
We first set up the scattering problem for an angle of incidence of $30^{\circ}$.
The original matrix is of dimension $7040 \times 7040$, which is compressed to
one of dimension $422 \times 452$.
The time for compression was 416.6 secs.\ (and for generating the
matrix entries, 1762.3 secs.).
The time for inversion was 2.9 secs.
The relative error in the solution (compared with standard LU factorization)
was $1.23\, \times 10^{-6}$. We then changed the angle of incidence to
$45^{\circ}$ and used the updating method of section \ref{multangle}.
The time for updating the compressed forward operator was 68.4 secs.\
and the relative error in the solution was $7.13\, \times 10^{-6}$.
In this problem, there are six propagating modes, with directions indicated
in Fig.~\ref{scatmodes2}.

\vspace{.1in}

\begin{figure}[!htb]   
\begin{center}
\includegraphics[width=3.5in]{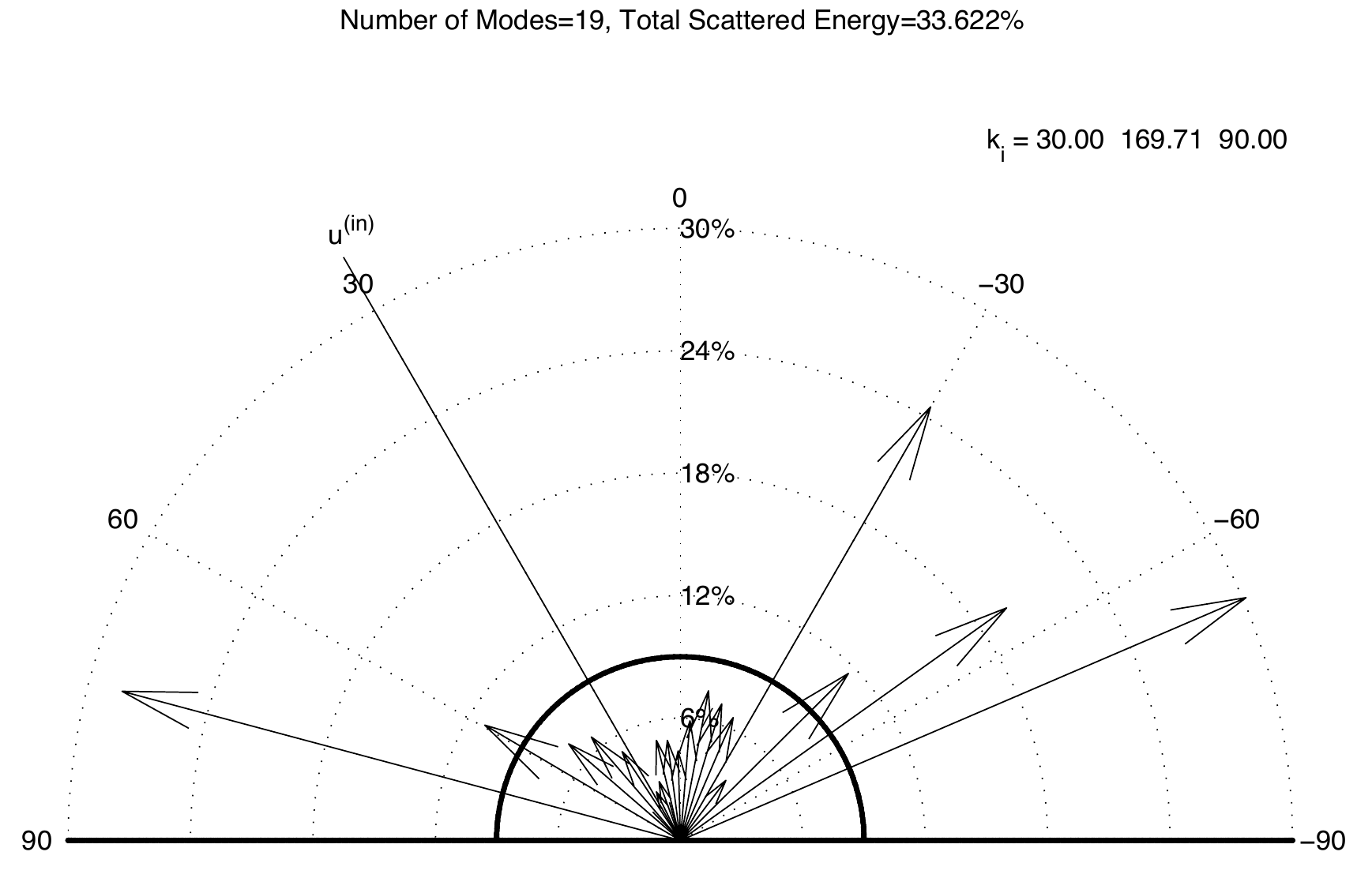}
\end{center}
\caption{The strength of the 19 radiation modes in Example 3 with a
semicircular scatterer with $k_0=30$, $k_1 = 120 \sqrt{2}$, $k_2 = 90$ and $d=2$.
Only the scatterer in the unit cell is plotted.
The length of the arrows in the various diffraction directions indicate the relative
amplitudes $a^+_n$. Approximately
$33.622 \%$ of the energy is scattered upward.
\label{hemifig}}
\end{figure}

\noindent
{\bf Example 3:}
The complexity of the scattering pattern can be quite striking.
In Fig.~\ref{hemifig} is shown the scattering pattern
from a semicircular scatterer with an angle
of incidence of $30^{\circ}$.
We set $\omega=10$ and choose $\epsilon$ so that the
Helmholtz coefficient in the upper half-space, the semicircular scatterer and the
substrate layer are $k = 30$, $120 \sqrt{2}$, and $90$,
respectively. There are 19 radiation modes at this angle of incidence.

\vspace{.1in}

\begin{figure}[!htb]   
\begin{center}
\includegraphics[width=5in]{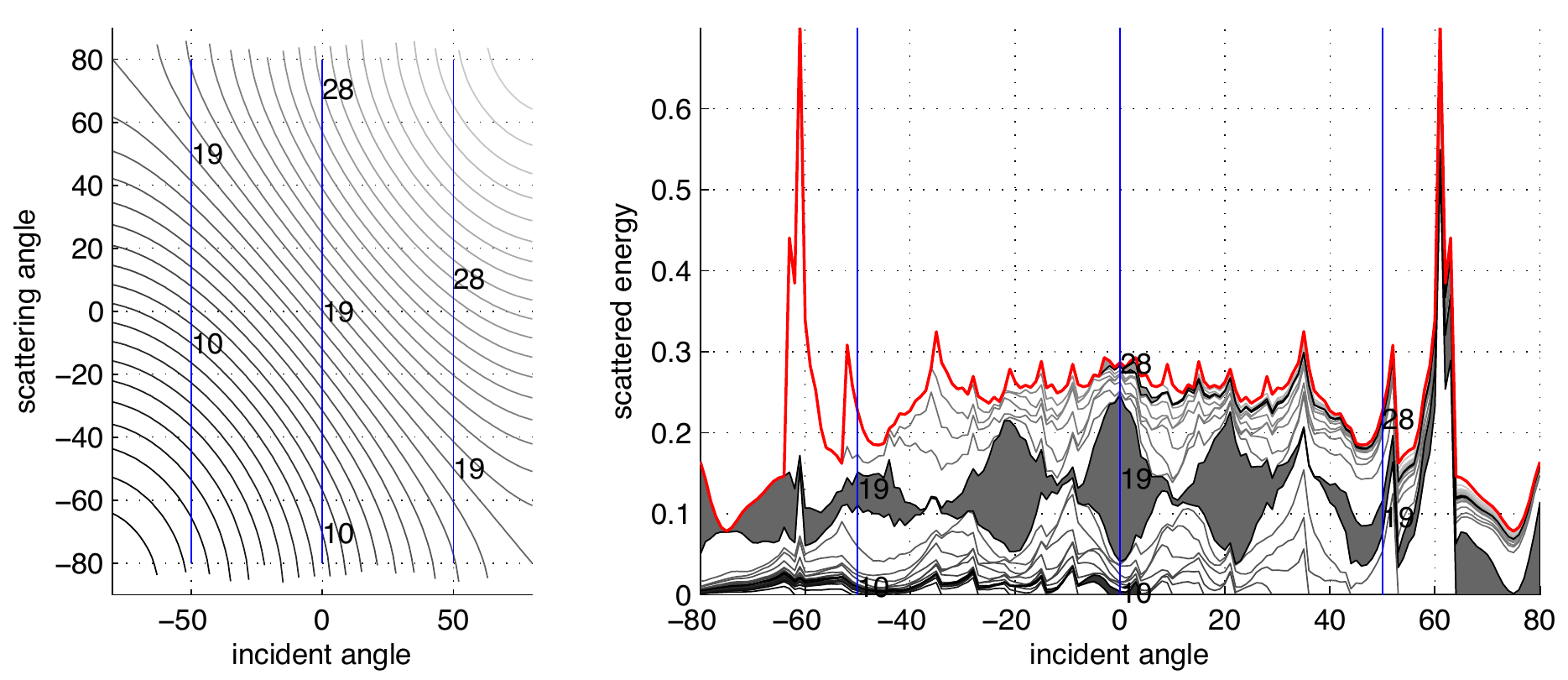}
\end{center}
\caption{The Bragg scattering angles as a function of incident angle
(left) and the scattering energies in the various modes for the
trapezoidal-shaped scatterer shown in Example 1, with $k_0=30$.
(See the text for a discussion of the plots.)
\label{trapscan}}
\end{figure}

\begin{figure}[!htb]   
\begin{center}
\includegraphics[width=5in]{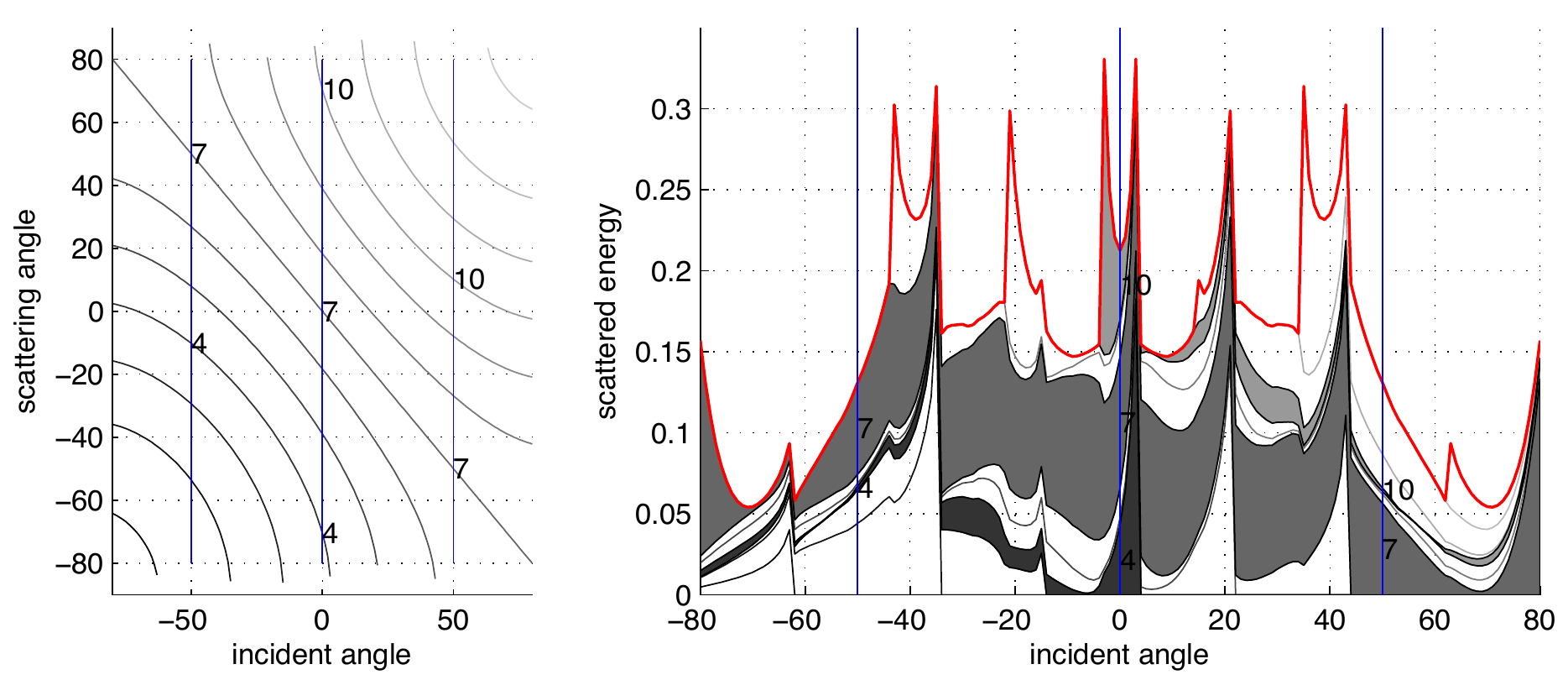}
\end{center}
\caption{The Bragg scattering angles as a function of incident angle
(left) and the scattering energies in the various modes for the
semicircular scatterer of Example 3, with $k_0=10$.
(See the text for a discussion of the plots.)
\label{hemiscan}}
\end{figure}

\noindent
{\bf Examples 4, 5:}
In our final examples, we compute the diffraction pattern across all angles
of incidence from $\theta = -80^{\circ}$ to $\theta = 80^{\circ}$ for the scattering geometries
depicted in Examples 1 and 3, except that for the trapezoidal-shaped scatterer, we increased
$\omega$ by a factor of 3, so that $k_0=30$ instead of 10.
For the semicircular scatterer, we decreased
$\omega$ by a factor of 3, so that $k_0=10$ instead of 30.
On the left-hand side of Figs.~\ref{trapscan} and \ref{hemiscan} are plotted the
diffraction orders as a function of incident angle. That is,
for each incident angle $\theta$, the intersection of the indicated vertical line
with the various curves are the {\em Bragg angles}
$\theta_n = \tan^{-1}(k_n/\kappa_n)$ according to formula
(\ref{e:rbu}), where $k_n$ and $\kappa_n$ are chosen to enforce both quasiperiodicity and
the Helmholtz equation.

\begin{remark}
The number of intersections of each vertical line on these left-hand plots
defines the precise number of modes for a given angle of incidence.
It is easy to see that each of the curves on the left-hand plots
traverses the {\em incident angle}-{\em scattered angle} plane
continuously (until it disappears), so that we may enumerate the modes unambiguously
from the lower left
corner to the upper right corner. The labels (``10", ``19", ``28") in
Fig.~\ref{trapscan} are drawn on the 10th, 19th, and 28th such curve.
The labels (``4", ``7", ``10") in
Fig.~\ref{hemiscan} are drawn on the 4th, 7th, and 10th such curve.
\end{remark}

On the right-hand side of Figs.~\ref{trapscan} and \ref{hemiscan} are plotted the
fraction of energy radiated into each mode. The $i$th curve in the right-hand plots
show the total energy scattered in modes $1$ through $i$. Thus, the separation
between curves corresponding to modes $i$ and $(i-1)$ shows the fraction of energy
radiated in the $i$th mode. Highlighted in gray are the energies scattered in the
10th, 19th, and 28th modes in
Fig.~\ref{trapscan} and in the
4th, 7th, and 10th modes in Fig.~\ref{hemiscan}.
Note that the strength can change quite abruptly when the incident angle is changed
only slightly.

\section{Conclusions} \label{sec:conclusions}

We have described an integral equation method for quasi-periodic scattering
from layered materials with grating-like structures on the ``top" surface.
It combines
(1) the use of the quasi-periodic Green's function,
(2) the modified Kress/M\"{u}ller/Rokhlin integral equation
for multi-material junctions \cite{triplepoint},
(3) the use of exponential refinement near geometric singularities,
\cite{BRS,HO}, and (4) the fast direct solver of \cite{HG}.

Since the quasi-periodic Green's function changes with each incident angle,
there is a global change to the system matrix with each new illumination.
We have shown, however, that
the {\em difference} between
Green's functions at different angles of incidence is (hierarchically) smooth so that
the compressed representation of the system matrix can be rapidly updated.

In recent work, Gillman and Barnett \cite{gillmanbarnett} developed an
alternative fast direct solver based on using the free-space Green's function
with auxilliary variables to impose quasi-periodicity. In that formulation,
the bulk of the matrix is left unchanged for different illuminations.
We suspect that the relative advantages of the two approaches will depend
on the aspect ratio of the unit cell, the spatial dimension (2D vs.\ 3D scattering)
and detailed implementation issues.
Both approaches have asymptotically optimal complexity
for unit cells that are a modest number of wavelengths in size.

\section*{Acknowledgements}
This work was supported by the Applied
Mathematical Sciences Program of the U.S. Department of Energy
under Contract DEFGO288ER25053 and
by the Air Force Office of Scientific Research under
NSSEFF Program Award FA9550-10-1-0180.
KLH was also supported in part by the National Science Foundation under grants DGE-0333389 and DMS-1203554.
JYL was also supported in part by the Priority Research Centers Program
(2009-0093827) and the Basic Science Research Program (2012-002298)
through the National Research Foundation (NRF) of Korea.

\bibliographystyle{plain}

\end{document}